\documentclass{article}
\usepackage[pdftex]{graphicx}   
\usepackage{amssymb}
\usepackage{color}

\newcommand{\ord}{\mathrm{ord}\,}
\newcommand{\cont}{\mathrm{cont}}
\newcommand{\lc}{\mathrm{lc}}
\newcommand{\ini}{\mathrm{in}\,}
\newcommand{\inN}{\mathrm{in}_{\cal N}}
\newcommand{\Zer}{\mathrm{Zer}\,}
\newcommand{\bC}{\mathbf{C}}

\newcommand{\bR}{\mathbf{R}}

\newcommand{\bU}{\mathbf{U}}
\newcommand{\cD}{{\cal D}}
\newcommand{\cC}{{\cal C}}
\newcommand{\classB}{\overline{B}}

\newtheorem{Theorem}{Theorem}[section]
\newtheorem{Lemma}[Theorem]{Lemma}
\newtheorem{Proposition}[Theorem]{Proposition}
\newtheorem{Corollary}[Theorem]{Corollary}
\newtheorem{Definition}[Theorem]{Definition}
\newtheorem{Remark}[Theorem]{Remark}
\newtheorem{Property}[Theorem]{Property}
\newtheorem{Formula}[Theorem]{Formula}

\newtheorem{Example}[Theorem]{Example}

\newenvironment{proof}[1][Proof]{\textbf{#1.} }{\
\rule{0.5em}{0.5em}}

\newcommand{\Teis}[2]{
   \setlength{\unitlength}{1ex}
   \begin{picture}(2,0)(0,0.4)
      \put(0,1.1){\line(1,0){2}}
      \put(0,0.9){\line(1,0){2}}
      \put(1,1.2){\makebox(0,0)[b]{$\scriptstyle #1$}}
      \put(1,0.8){\makebox(0,0)[t]{$\scriptstyle #2$}}
   \end{picture}}

\newcommand{\Teisssr}[4]{
   \setlength{\unitlength}{1ex}
   \begin{picture}(#3,3)(0,0.4)
      \put(0,1.15){\line(1,0){#3}}
      \put(0,0.85){\line(1,0){#3}}
      \put(#4,1.3){\makebox(0,0)[b]{$#1$}}
      \put(#4,0.7){\makebox(0,0)[t]{$#2$}}
   \end{picture}}

\title{Non-degeneracy of the discriminant  \footnotetext{
     \noindent   \begin{minipage}[t]{5in}
       {\small
       2000 {\it Mathematics Subject Classification:\/} Primary 32S55;
       Secondary 14H20.\\
       Key words and phrases: irreducible plane curve,
       jacobian Newton diagram, polar invariant, discriminant, degeneracy.\\
       The first-named author was partially supported by the Spanish
       Project PNMTM 2007-64007. The first-named and third-named authors were
       partially supported by the Polish MSHE grant No N N201 386634.}
       \end{minipage}}}

\author{Evelia R.\ Garc\'{\i}a Barroso, Janusz Gwo\'zdziewicz and
Andrzej Lenarcik}

\begin{document}
\maketitle
\begin{center}
{\it To Professor Arkadiusz P\l oski on his
65$^{\mbox{\scriptsize th}}$ birthday}
\end{center}
\begin{abstract}
\noindent  Let $(\ell ,f)\colon (\mathbf C^{2},0)\longrightarrow (\mathbf
C^2,0)$ be the germ of a holomorphic mapping such that  $\ell=0$ is a 
smooth curve and $f=0$ has an isolated singularity at $0\in \bC^2$.  We assume that
$\ell=0$ is not a branch of $f=0$.
The direct image
of the critical locus of this mapping is called the discriminant curve.   In this paper 
we study the pairs $(\ell,f)$  for which the discriminant curve is non-degene\-rate in  the
Kouchnirenko sense.
\end{abstract}

\section{Introduction}

\medskip

\noindent Let $(\ell ,f)\colon (\mathbf C^{2},0)\longrightarrow (\mathbf
C^2,0)$ be a holomorphic mapping given by $u=\ell (x,y)$, $v=f(x,y)$, where $\ell=0$ is a
smooth curve and $f=0$ has an isolated singularity at $0\in \bC^2$.  We assume that
$\ell=0$ is not a branch of $f=0$.
To any such morphism we can associate two analytic curves: the {\em polar curve} 
$\frac{\partial \ell}{\partial x}\frac{\partial f}{\partial
y} - \frac{\partial \ell}{\partial y}\frac{\partial f}{\partial x}=0$ and its direct image $\cD(u,v)=0$ 
which is called the {\em discriminant curve} of the morphism $(\ell ,f)$  (see
\cite{Teissier1}, \cite{Casas-Asian}). A series $\cD(u,v)$, defined up to multiplication by an invertible power series, 
is called the {\em discriminant}. In \cite{Teissier1} and
\cite{Teissier3} Teissier introduced the {\em jacobian Newton
diagram}, which is the Newton diagram of $\cD(u,v)$. The jacobian Newton diagram depends only on the
topological type of $(\ell,f)$ (see \cite{Teissier1} for the case where $\ell$ is generic, Merle \cite{Merle} and Ephraim \cite{Ephraim} for one branch and \cite{Gwo1},
\cite{Lenarcik} and \cite{Michel} for general case). Decompositions of the polar curve
can be found in the literature (see \cite{Merle}, \cite{Ephraim},  \cite{Eggers}, \cite{GB}) . In the spirit of  Eggers  \cite{Eggers} we propose a factorization of  the discriminant $\cD (u,v)$. The Newton diagram of every factor has only one compact edge.
We specify formulas for the weighted initial forms of these factors. Using this description we study the  pairs $(\ell,f)$  for which the discriminant is non-degene\-rate, in the  Kouchnirenko sense
\cite{Kouchnirenko}, answering a question of Patrick Popescu-Pampu.

\medskip

\noindent For the irreducible case we prove in Section \ref{irreducible}:

\begin{Theorem}
\label{teorema-ramas} Let $f=0$ be a branch. Then the discriminant of $(\ell,f)$  is non-degenerate
if and only if there are no lattice
points inside the compact edges of its Newton diagram.
\end{Theorem}

\begin{Corollary}
\label{topological type}  Let $f=0$ be a branch. Then the non-degeneracy of the discriminant of $(\ell,f)$  depends only on
the topological type of $(\ell,f)$.
\end{Corollary}

\noindent Theorem~\ref{teorema-ramas} is not true for reducible curves. 
In Examples~\ref{ex-1} and \ref{ex0} we construct two equisingular  (that is having the same embedded topological type) curves $f_1=0$ and $f_2=0$ such that the 
jacobian Newton diagrams of morphisms $(x,f_1)$ and $(x,f_2)$ have only one compact edge 
joining points $(0,4)$ and $(40,0)$. This edge has three lattice points inside. Nevertheless the discriminant 
of $(x,f_1)$ is degenerate while that of $(x,f_2)$ is non-degenerate. 

\medskip
\noindent The structure of the paper is as follows: in Section \ref{Preliminaries}  we start  by recalling the notion of non-degeneracy. Then, after a change of coordinates, we may assume 
that the morphism that we consider has the form $(x,f)$. We describe the discriminant
 by using Newton--Puiseux roots of the $y$-partial derivative of $f(x,y)$. For that the Lemma of Kuo-Lu  plays an important role. Using the results of this section we construct  examples of curves with many smooth branches, which determine non-degenerate discriminants.

\medskip
\noindent In Section \ref{section-factorization} we propose an analytical factorization of $\cD (u,v).$ In Proposition \ref{initial-disc} we compute the initial Newton polynomial of every factor and express it as a product of rational powers of quasi-homogeneous polynomials. Then in Section \ref{irreducible} we apply this formula to irreducible power series $f(x,y)$  and we characterize in Corollary \ref{caract-branches}
the equisingularity  classes of branches for which  the discriminant of $(x,f)$ is non-degenerate.

\medskip
\noindent In Section \ref{general-section} we return to the general case.~Taking up again Proposition \ref{initial-disc} we give, in Proposition \ref{polynomial-factorization}, a polynomial factorization of the initial Newton polynomials of the factors of  $\cD (u,v)$. As a consequence, in Theorem \ref{Th:general},  we obtain a criterion for non-degeneracy of the factors of the discriminant.
We finish this section with another example of curves with as many singular branches as we wish, which determine non-degenerate discriminants. 

\medskip
\noindent In the last section we analyze what impact on the discriminant has a modification of $\ell$ or $f$ in the morphism $(\ell,f)$. 
Theorem~\ref{tc1} shows that non-degeneracy of the discriminant of the morphism 
$(\ell,f)$ is independent of the choice of the representative  of the curve $f=0$. 
Theorem~\ref{tc3} shows that if $f=0$ is unitangent and transverse to $\ell=0$, then the
non-degeneracy of the discriminant of the morphism $(\ell,f)$ depends only on the curve $f=0$.
The assumption that $f=0$ has only one tangent cannot be omitted as it  is shown in Example~\ref{ejemJ}.

\section{Preliminaries}
\label{Preliminaries}

\noindent We start this section recalling the notion of non-degeneracy. Then we reduce our study to the morphisms of the form $(x,f)$. We describe the discriminant
 by using Newton-Puiseux roots of $\frac{\partial f}{\partial y}(x,y)$. 
The Lemma of Kuo-Lu  plays an important role.

\subsection{Non-degeneracy after Kouchnirenko}
\label{Newton polygons}
\medskip

\noindent 
Set $\bR_{+}=\{\,x\in\bR: x\geq0\,\}$. Let 
$f(x,y)=\sum_{ij} a_{ij}x^iy^j\in\bC\{x,y\}\backslash\{0\}$.
The {\em Newton diagram} of $f$ is
$$
\Delta_f:= \mbox{Convex Hull}\;\left(\{\,(i,j): a_{ij} \neq0\,\}+\bR_{+}^2\right).
$$

\noindent The
Newton diagram of a product is the Minkowski sum of the Newton diagrams of the factors.
That is $\Delta_{fg}=\Delta_f+\Delta_g$,  where
$\Delta_f+\Delta_g=\{\,a+b:a\in\Delta_f, b\in\Delta_g\,\}$. In
particular if $f$ and $g$ differ by an invertible factor $u\in
\mathbf C\{x,y\}$, $u(0,0)\neq0$ then $\Delta_f=\Delta_g$. 

\medskip

\noindent The {\em initial Newton polynomial} of $f(x,y)=\sum_{i,j} a_{ij}x^iy^j$, denoted by $\inN f$, is the sum of all terms $a_{ij}x^iy^j$ such that $(i,j)$ belongs to a compact edge of $\Delta_f$.

\medskip

\noindent Following Teissier \cite{Teissier2} we introduce {\em
elementary Newton diagrams}. For $m,n>0$  we put
$\{\Teis{n}{m}\}=\Delta_{x^n+y^m}$. We put also
$\{\Teis{n}{\infty}\}=\Delta_{x^n}$ and
$\{\Teis{\infty}{m}\}=\Delta_{y^m}$. By definition the {\em
inclination} of $\{\Teis{L}{M}\}$  is $L/M$ with the conventions that
$L/\infty=0$ and $\infty/M=+\infty$.

\medskip

\noindent Let $S$ be a compact edge of ${\Delta}_f$ of inclination
$p/q$,  where $p$ and $q$ are coprime integers. The {\em
initial part of $f(x,y)$ with respect to $S$} is the
quasi-homogeneous polynomial $f_S(x,y)=\sum
a_{ij}x^i y^j$ where the sum runs over all lattice points  $(i,j)\in S$. Observe that if $\Delta_f$ is an elementary Newton diagram then the initial part of $f(x,y)$ with respect to the only compact edge of $\Delta_f$ coincides with the initial Newton polynomial of $f(x,y)$.

\medskip

\noindent Decomposing $f_S(x,y)$ into irreducible factors in
$\mathbf C[x,y]$ we get
\begin{equation}
\label{factorization} f_S(x,y)=cx^k y^l
\prod_{i=1}^r(y^q-a_ix^p)^{s_i},
\end{equation}

\noindent where  $k$ and $l$ are non-negative integers, $c$ and $a_i$ are nonzero complex numbers and $a_i\neq a_j$ for
$i\neq j$.

\medskip

\noindent The series
$f(x,y)$ is {\em non-degenerate on the compact edge} $S$ of ${\Delta}_f$ if in  (\ref{factorization}) $s_i=1$ for all
$i\in \{1,\ldots,r\}$. In particular $f$ is non-degenerate on the compact edge $S$ if there are no lattice points inside $S$. The converse is not true
as  $(y-x)(y-2x)$ shows. The series $f(x,y)$
is  {\em non-degenerate} if it is non-degenerate on every compact edge of its Newton diagram (see \cite{Kouchnirenko}).

\medskip

\subsection{Newton--Puiseux roots}
\label{Notations}
\medskip

\noindent  Let $\bC\{x\}^*$ be the ring of Puiseux series in $x$, that is
the set of  series of the form
$$ \alpha(x)=a_1x^{N_1/D}+a_2x^{N_2/D}+\cdots,
   \quad a_i\in\bC,
$$

\noindent where $N_1<N_2<\dots$ are non-negative integers, $D$ is a positive integer and
$a_1t^{N_1}+a_2t^{N_2}+\cdots$ has a positive radius of convergence.
In this paper $+\cdots$ means {\em plus higher order terms}. If
$a_1\neq0$ then  the {\em order} of $\alpha(x)$ is
$\ord\alpha(x)=N_1/D$ and the {\em initial part} of $\alpha(x)$ equals $\ini \alpha(x)=a_1x^{N_1/D}$. By convention the order of the zero series is
$+\infty$.
For any Puiseux series $\alpha(x)$, $\gamma(x)$ we denote by
$O(\alpha,\gamma)=\ord(\alpha(x)-\gamma(x))$ and call this number
the {\em contact order} between $\alpha(x)$ and $\gamma(x)$. If $Z\subset\bC\{x\}^*$ is a finite
set then the \textit{contact} between
$\alpha\in\bC\{x\}^*$ and $Z$ is $\cont(\alpha,Z)=\max_{\gamma\in Z}
O(\alpha,\gamma)$.

\medskip

\noindent By a {\em fractional power series} we mean a Puiseux series of positive order.

\medskip

\noindent Let $g(x,y)\in\bC\{x,y\}$ be a convergent power series. A
fractional power series  $\gamma(x)$ is called a
{\em Newton--Puiseux root} of $g(x,y)$ if $g(x,\gamma(x))=0$ in
$\mathbf C\{x\}^*$. We denote by $\Zer g$ the set of all Newton--Puiseux
roots of $g(x,y)$.

\medskip
\noindent If $g=g_1^{a_1}\cdots g_r^{a_r}$ where the $g_i$ are irreducible and pairwise coprime  elements of $\bC\{x,y\}$, then the curves
$g_i=0$ are called the {\em branches} of $g=0$. We
say that $g=0$ is {\em reduced} if $a_1=\cdots=a_r=1$. Notice that $g$ has an isolated singularity at $0\in \bC^2$ if and only if it is singular and reduced.

\subsection{The Lemma of Kuo-Lu}

\noindent Consider the morphism $(\ell,f)$ as in  Introduction, where $f$ is a reduced curve. An analytic change of coordinates does not affect the discriminant curve (see for example \cite{Casas-Asian}, Section 3). Hence in what follows  we assume that $\ell(x,y)=x$. Then $\frac{\partial f}{\partial y}=0$ is the polar curve of $(x,f)$. 
\medskip

\noindent
The Newton--Puiseux factorizations of $f(x,y)$ and $\frac{\partial f}{\partial y}(x,y)$ are of the form
\begin{equation}
\label{ppp1}
f(x,y)= u(x,y)\prod_{i=1}^p [y-\alpha_i(x)], 
\end{equation}
\begin{equation}
  \label{ppp2}
\frac{\partial f}{\partial y}(x,y)=\tilde u(x,y)\prod_{j=1}^{p-1} [y-\gamma_j(x)],
\end{equation}
where 
$u(x,y)$, $\tilde u(x,y)$ are units in  $\bC\{x,y\}$ and $\alpha_i(x)$, $\gamma_j(x)$
are fractional power series. Since $f$ is reduced, 
$\alpha_i(x)\neq \alpha_j(x)$ for $i\neq j$.

\medskip

\noindent The following lemma, which is a part of Lemma 3.3 in
\cite{Kuo-Lu} (for the transverse case; see \cite{G}, Corollary 3.5  and \cite{G-P}, Proposition 2.2 for the general
case), describes the contacts between Newton--Puiseux roots of
$f(x,y)$ and $\frac{\partial f}{\partial y}(x,y)$.

\begin{Lemma}\label{Wn2}
For every $\gamma_j\in \Zer\frac{\partial f}{\partial y} $ there exist $\alpha_k,\alpha_l\in\Zer f$, $k\neq l$ such that
$$ O(\alpha_k,\gamma_j) = O(\alpha_l,\gamma_j) =  O(\alpha_k,\alpha_l)=
     \max_{i=1}^p\, O(\alpha_i,\gamma_j) .
$$
\end{Lemma}

\noindent In what follows we recall the {\em tree model} introduced in \cite{Kuo-Lu} which encodes the contact orders
between Newton--Puiseux roots of $f(x,y)$.

\begin{Definition}
Let $\alpha\in\bC\{x\}^*$  and
let $h$ be a positive rational number.
The pseudo-ball $B(\alpha,h)$ is the set
$ B(\alpha,h)=\{\,\gamma\in\bC\{x\}^*: O(\gamma,\alpha)\geq h \,\} .
$
We call $h(B):=h$ the height of $B:=B(\alpha,h)$.
\end{Definition}

\medskip

\noindent Note that $h(B)$ is well-defined since $h(B)=\inf\{O(\gamma,\beta)\;:\;\gamma,\beta \in B\}$.

\medskip

\noindent Consider the following set of pseudo-balls
$$T(f):=\{ B(\alpha, O(\alpha,\alpha'))\;:\; \alpha, \alpha'\in\Zer f,\,
\alpha\neq \alpha' \}.$$

\noindent The elements of $T(f)$ can be identified with bars of the
{\em tree model} of~$f$ defined in \cite{Kuo-Lu} (for a short  presentation see also Section 8 of \cite{I-K-K}). It follows from
Lemma~\ref{Wn2} that for every $\gamma \in \Zer\frac{\partial
f}{\partial y}$ there exists exactly one $B\in T(f)$ such that
$\gamma\in B$ and $h(B)=\cont(\gamma, \Zer f)$. Following
\cite{Kuo-Parusinski2} we say that $\gamma$ {\em leaves} $T(f)$ at $B$.

\medskip

\noindent Take a pseudo-ball $B \in T(f)$ .
Every $\gamma\in B$ has the form

\begin{equation}
\label{expression}
 \gamma(x)=\lambda_B(x)+c_{\gamma}x^{h(B)}+\cdots,
\end{equation}

\noindent 
where $\lambda_B(x)$  is obtained from an arbitrary $\alpha(x)\in B$ by omitting all the terms of order bigger than or equal to $h(B)$.

\medskip

\noindent We call the
 complex number $c_{\gamma}$  the {\em leading
coefficient} of $\gamma$ with respect to $B$ and we
denote it by $\lc_B(\gamma)$.  Remark that   $c_{\gamma}$ can be zero.

\medskip

\noindent We need next two Lemmas from \cite{G}  (see also \cite {LMP} and \cite{Lenarcik-2004}, Corollary 3.7 and Proposition 3.6).

\begin{Lemma}[\cite{G}, Lemma 3.3]
 \label{L2}
Let $B \in T(f)$.
There exist a polynomial $F_B(z)\in\bC[z]$, depending on $f$, and a rational number $q(B)$
such that for every $\gamma(x)=\lambda_B(x)+c_{\gamma}x^{h(B)}+\cdots$
\begin{equation}\label{Eq4}
f(x,\gamma(x))=F_B(c_{\gamma})x^{q(B)}+\cdots\; .
\end{equation}
Moreover
\begin{equation}\label{Eq10}
F_B(z)=C\!\!\prod_{i:\alpha_i\in B} (z-\lc_B(\alpha_i)),
\end{equation}where $C$ is a nonzero constant.
\end{Lemma}

\begin{Remark}
\label{xxx}
It follows from the proof of Lemma 3.3 in \cite{G} that if $f$ is a Weierstrass polynomial and $\alpha_j(x)\in B$, then the constant $C$ in (\ref{Eq10}) is expressed by the formula
$$Cx^{q(B)}=\prod_{i:\alpha_i\not\in B}\ini (\alpha_j(x)-\alpha_i(x))\prod_{i:\alpha_i \in B}x^{h(B)}.$$
\end{Remark}
\begin{Lemma}[\cite{G}, Lemma 3.4]

\label{L1}  Let $B \in T(f)$.  Then
$$
\frac{d}{dz} F_B(z)=
C' \prod_{j:\gamma_j\in B} (z-\lc_B(\gamma_j)),
$$

\noindent where $C'$ is a nonzero constant.
\end{Lemma}

\medskip

\medskip

\noindent Using the above lemmas we characterize the  Newton--Puiseux roots
of $\frac{\partial f}{\partial y}(x,y)$  leaving $T(f)$ at a fixed $B$.

\begin{Lemma}
\label{caract}
Let $B\in T(f)$ and $\gamma \in B$. Then $\gamma$ leaves $T(f)$ at $B$ if and only if
$F_B(\lc_B(\gamma))\neq 0$.
\end{Lemma}
\noindent \begin{proof}
For $\gamma \in B$ the inequality $F_B(\lc_B(\gamma))\neq 0$ is
equivalent to $\lc_B(\gamma)\neq\lc_B(\alpha_i)$
for all $\alpha_i\in B$, and this is equivalent to $\cont(\gamma,\Zer f)=h(B)$.
\end{proof}

\medskip 

\noindent Given $B,B'\in T(f)$, we say that $B'$ is a  {\em direct successor} of $B$ in $T(f)$ if
$B\supset B'$ and  there are no $B''\in T(f)$ (different from $B$ and $B'$) such that $B\supset B'' \supset B'$. 
The next lemma follows from  Theorem C in \cite{Kuo-Lu}. For convenience of the reader we present a proof:  

\begin{Lemma}
\label{successor}
Let $B,B'\in T(f)$. Suppose that $B'$ is a direct successor of $B$ in $T(f)$. Then 
$q(B')-q(B)=\sharp (B'\cap \Zer f)[h(B')-h(B)]$, where the symbol $\sharp$ stands for the number of the elements of a set. If $B\in T(f)$ is the pseudo-ball of the minimal height then $q(B)=\sharp (\Zer f)h(B).$
\end{Lemma}
\noindent \begin{proof}
Let $\delta(x)=\lambda_B(x)+cx^{h(B)}$ where $F_B(c)\neq 0$ and 
$\delta'(x)=\lambda_{B'}(x)+c'x^{h(B')}$ where $F_{B'}(c')\neq 0$. Then following  (\ref{ppp1}) and
Lemma \ref{L2}
\begin{equation}
\label{aaa1}
q(B)=\ord f(x,\delta(x))=\sum_{\alpha\in \Zer f}O(\delta,\alpha)
\end{equation}
and 
\begin{equation}
\label{aaa2}
q(B')=\ord f(x,\delta'(x))=\sum_{\alpha\in \Zer f}O(\delta',\alpha)
.\end{equation}

\noindent  We have $O(\delta,\alpha)=h(B)$, 
$O(\delta',\alpha)=h(B')$ for $\alpha \in \Zer f \cap B'$.   
U\-sing the strong triangle inequality property of the contact order
one checks that $O(\delta,\alpha)=O(\delta',\alpha)$ for  $\alpha \in \Zer f\backslash B'$.
Substracting (\ref{aaa1}) from (\ref{aaa2})  we get the first statement of the lemma. 
The second statement of the lemma is  a consequence of (\ref{aaa1}).
\end{proof}

\medskip

\noindent Following Lemma 5.4 in \cite{GB-G} the discriminant of
the morphism $(x,f)$ can be written as
\begin{equation}
\label{equation-discriminant}
\cD(u,v) = \prod_{j=1}^{p-1}(v-f(u,\gamma_j(u))) .
\end{equation}

\begin{Example}
\label{ex-1}
{\rm  Let $h(x,y)=(y-x^2-x^3)(y-x^2+x^3)(y+x^2-x^3)(y+x^2+x^3)$ 
and let $f_1(x,y)=x^{10}+\int_0^y h(x,t)\, dt=0$. 
Since $\frac{\partial f_1}{\partial y}(x,y)=h(x,y)$, we get by~(\ref{equation-discriminant})
$\inN\cD(u,v)=\bigl(v-\frac{23}{15}u^{10}\bigr)^2\bigl(v-\frac{7}{15}u^{10}\bigr)^2$. 
Thus the discriminant of $(x,f_1)$ is degenerate. 
One can also show that it remains degenerate after any analytical change of coordinates.}
\end{Example}

\begin{Example}
\label{ex0}
{\rm  Let $f_2(x,y)=y^5+x^8y+x^{10}$. As $f_2(x,y)$ is a quasi-homogeneous polynomial, all its 
Newton--Puiseux roots are monomials of the same order. The same applies to $\frac{\partial f_2}{\partial y}$. 
The tree model $T(f_2)$ has only one pseudo-ball $B$ of the height~2. We have 
$F_B(z)=f_2(1,z)=z^5+z+1$. All critical values $w_j=F_B(z_j)$, 
where $z_1$, \dots , $z_4$ are critical points of $F_B(z)$, are pairwise different. 
By~(\ref{equation-discriminant}) and Lemma~\ref{L1} we get 
$\cD(u,v) =\prod_{j=1}^{4}(v-f_2(u,z_ju^2)))=
\prod_{j=1}^{4}(v-w_ju^{10})$. Hence the discriminant of $(x,f_2)$ is non-degenerate.}
\end{Example}

\noindent The curves $f_1(x,y)=0$ and $f_2(x,y)=0$ are equisingular. Nevertheless the discriminant of $(x,f_1)$ is degenerate while the discriminant of $(x,f_2)$ is non-degenerate.

\begin{Example}
\label{ex1}
{\rm \noindent Let $f(x,y)=\prod_{i=1}^4(y-\alpha_i(x))$ where  
$\alpha_1(x)=x+x^3$, $\alpha_2(x)=x-x^3$,
$\alpha_3(x)=-x+x^4$ and $\alpha_4(x)=-x-x^4$. The curve $f=0$ has four smooth branches.

\medskip
\noindent  The tree model $T(f)$ is given in the picture below.
Following \cite{Kuo-Lu} we draw pseudo-balls of finite height as horizontal bars.
The tree $T(f)$ has three bars: $B_1$ of height $1$, $B_2$ of height~$3$ 
and $B_3$ of height~$4$.

$$
\begin{picture}(110,65)(0,0)
\put(55,0){\line(0,1){15}}  
{\thicklines \put(15,15){\line(1,0){75}}}
\put(93,13){$B_1$}
{\thicklines \put(0,35){\line(1,0){30}}}
\put(15,15){\line(0,1){20}}  
\put(0,35){\line(0,1){20}} 
\put(30,35){\line(0,1){20}}
\put(0,55){$\alpha_1$} 
\put(30,55){$\alpha_2$}
\put(33,33){$B_2$}
{\thicklines \put(75,40){\line(1,0){30}}}
\put(90,15){\line(0,1){25}} 
\put(75,40){\line(0,1){20}} 
\put(105,40){\line(0,1){20}}
\put(75,60){$\alpha_3$} 
\put(105,60){$\alpha_4$}
\put(108,38){$B_3$}
\end{picture}
$$

\medskip

\noindent In order to compute the polynomial $F_{B}(z)$ for $B\in T(f)$ it is enough to
find the lowest  order term of $f(x,\lambda_{B}(x)+zx^{h(B)})$.

\medskip 

\noindent Since $\lambda_{B_1}(x)=0$ and $h(B_1)=1$, we get $f(x,\lambda_{B_1}(x)+zx^{h(B_1)})=f(x,zx)=(z-1)^2(z+1)^2x^4+{}\cdots$.

\medskip

\noindent  Similarly $f(x, \lambda_{B_2}(x)+zx^{h(B_2)})=f(x,x+zx^3)=4(z-1)(z+1)x^8+{}\cdots$ and $f(x,\lambda_{B_3}(x)+zx^{h(B_3)})=
f(x,-x+zx^4)=4(z-1)(z+1)x^{10}+{}\cdots$.

\medskip

\noindent Hence 
$$\begin{array}{ll}
F_{B_1}(z)=(z-1)^2(z+1)^2,\quad                &q(B_1)=4,\\
F_{B_2}(z)=4(z-1)(z+1), \quad  &q(B_2)=8,\\
F_{B_3}(z)=4(z-1)(z+1),             &q(B_3)=10.
\end{array}$$

\noindent Each of the above polynomials has exactly two roots. Thus for every
$i\in\{1,2,3\}$
there exists  a unique critical point $z_i$, $F'_{B_i}(z_i)=0$ such that the critical value $w_i=F_{B_i}(z_i)$ is nonzero.
It follows from Lemmas~\ref{L1} and~\ref{caract} that
$z_i =\lc_{B_i}\gamma_i$  for some
$\gamma_i\in\Zer\frac{\partial f}{\partial y}$ which leaves $T(f)$ at $B_i$. 
By Lemma \ref{L2} we have $f(x,\gamma_i(x))=w_ix^{q(B_i)}+{}\cdots$.
In view of equality (\ref{equation-discriminant}) the initial Newton polynomial
of the discriminant $\cD(u,v)$ is the initial Newton polynomial of   $\prod_{i=1}^3(v-w_iu^{q(B_i)})$. 
Since this polynomial does not have multiple factors, the discriminant $\cD(u,v)$ is  non-degenerate. }
\end{Example}
\medskip

\noindent What matters in Example \ref{ex1} is that different $B\in T(f)$ have different $q(B)$
and also that $T(f)$ is a binary tree, 
hence for every $B\in T(f)$ the polynomial $F_B(z)$ has exactly two roots 
and consequently there exists exactly one 
$\gamma\in\Zer\frac{\partial f}{\partial y}$ which leaves $T(f)$ at $B$.  We use this idea in the next example.

\medskip
\begin{Example} 
\label{ex2}{\rm Let $g(x,y)$ be a power series which tree model $T(g)$ is presented in the figure below.
The numbers attached to the bars are the heights of corresponding pseudo-balls. 
Applying Lemma \ref{successor} one can check that $\{\,q(B):B\in T(g)\,\}=\{\,8,16,20,36,38,42,44\,\}$. 
By the same argument as before the discriminant of the morphism $(x,g)$ is non-degenerate.

$$
\begin{picture}(150,60)(0,0)
\put(75,0){\line(0,1){10}}  
{\thicklines \put(30,10){\line(1,0){90}}}
\put(122,8){1}
\put(30,10){\line(0,1){15}}  
{\thicklines \put(10,25){\line(1,0){40}}}
\put(52,23){3}
\put(10,25){\line(0,1){15}} 
{\thicklines \put(0,40){\line(1,0){20}}}
\put(22,38){13}
\put(0,40){\line(0,1){15}} 
\put(20,40){\line(0,1){15}}
\put(50,25){\line(0,1){15}} 
{\thicklines \put(40,40){\line(1,0){20}}}
\put(62,38){14}
\put(40,40){\line(0,1){15}} 
\put(60,40){\line(0,1){15}}
\put(120,10){\line(0,1){15}}  
{\thicklines \put(100,25){\line(1,0){40}}}
\put(142,23){4}
\put(100,25){\line(0,1){15}} 
{\thicklines \put(90,40){\line(1,0){20}}}
\put(112,38){15}
\put(90,40){\line(0,1){15}} 
\put(110,40){\line(0,1){15}}
\put(140,25){\line(0,1){15}} 
{\thicklines \put(130,40){\line(1,0){20}}}
\put(152,38){16}
\put(130,40){\line(0,1){15}} 
\put(150,40){\line(0,1){15}}
\end{picture}
$$

\noindent The curve $g=0$ from the above example decomposes into eight smooth branches. Following the idea of Example \ref{ex2} one can construct new examples of multibranched curves, with more levels in their tree models, whose  discriminants are non-degenerate.}
\end{Example}

\section{Factorization of the discriminant}
\label{section-factorization}
\noindent Assume that all the Newton--Puiseux roots of $f(x,y)$
and $\frac{\partial f}{\partial y}(x,y)$
belong to $\bC\{x^{1/D}\}$ for some positive integer $D$.
We define the action of the multiplicative group $\bU_D=\{\theta\in\bC\;:\;\theta^D=1\}$
of $D$-th complex roots of unity on $\bC\{x^{1/D}\}$.

\medskip

\noindent Take ~$\theta\in \bU_D$    and $\phi\in\bC\{x^{1/D}\}$ of the form
$\phi(x)=a_1x^{N_1/D}+a_2x^{N_2/D}+\cdots, $ where $0\leq N_1<N_2<\cdots $.
 By definition
$\theta \ast \phi(x)=a_1\theta^{N_1}x^{N_1/D}+a_2\theta^{N_2}x^{N_2/D}+\cdots$.
Following \cite{Kuo-Parusinski2} we call the series $\theta \ast \phi$ {\em conjugate} to $\phi$.

\medskip

\noindent It is well-known (see for example \cite{Walker}) that if $g(x,y)$ is an irreducible power series such that
$\Zer g\subset\bC\{x^{1/D}\}$ then the conjugate action of $\bU_D$ permutes transitively the
Newton--Puiseux roots of $g(x,y)$.
The conjugate action of $\bU_D$  pre\-serves  the contact order, i.e.
$O(\phi,\psi)=O(\theta\ast\phi,\theta\ast\psi)$ for $\phi$, $\psi\in\bC\{x^{1/D}\}$ and
$\theta\in \bU_D$.  

\medskip

\noindent The {\em index} of a fractional power series $\beta(x)$ is the smallest positive integer $N$ such that $\beta(x)\in \bC\{x^{1/N}\}$.  Following \cite{Walker} we get:

\begin{Property}
\label{index}
Let $\beta(x)\in \bC \{x^{1/D}\}$ be a fractional power series. Then the fo\-llo\-wing conditions are equivalent:
\begin{enumerate}
\item The index of $\beta(x)$ equals $N$.
\item The set $\{\theta \ast \beta(x)\;:\;\theta^D=1 \}$ has $N$ elements.
\item If $g(x,y)$ is an irreducible power series such that $g(x,\beta(x))=0$ then $\ord g(0,y)=N$.
\end{enumerate}
\end{Property}

\medskip

\noindent The action of $\bU_D$ on $\Zer f$ induces an action of this group on $T(f)$ as follows. Let 
$B=B(\alpha_k,{O}(\alpha_k,\alpha_l))$ and let $\theta\in\bU_D$.
Set $\theta\ast B=B(\theta\ast \alpha_k,{O}(\alpha_k,\alpha_l))$.
The properties of the conjugate action  imply that  $\theta\ast B$ is an element of $T(f)$ and
$\theta\ast B=B(\theta\ast\lambda_B,h(B))$.  Hence the definition of $\theta\ast B$ does not
depend on the choice of $\alpha_k\in B\cap \Zer f$.

\medskip

\begin{Proposition}
\label{invariant} Let $B\in T(f)$, $\theta\in\bU _D$  and $B'=\theta\ast B$. Then
$q(B)=q(B')$
and
$\theta^{q(B)D}F_B(z)=F_{B'}(\theta^{h(B)D}z)$.
\end{Proposition}

\noindent \begin{proof}
Acting by $\theta$ on the equation
$f(x,\lambda_B(x)+cx^{h(B)})=F_B(c)x^{q(B)}+\cdots$
we get
$f(x,\lambda_{B'}(x)+c\theta^{h(B)D}x^{h(B)})=
     F_B(c)\theta^{q(B)D} x^{q(B)}+\cdots\;.$
By Lemma~\ref{L2}
$f(x,\lambda_{B'}(x)+c\theta^{h(B)D}x^{h(B)})=
     F_{B'}(c\theta^{h(B)D}) x^{q(B')}+\cdots\;.$
\noindent Since $c$ is arbitrary, equating the right hand sides of the formulas above gives the proof.
\end{proof}

\bigskip

\noindent For every $B\in T(f)$ we denote by $\classB$ the orbit $\bU_D\ast B$ and by $E(f)$ the set of all orbits
in $T(f)$.

\medskip

\noindent 
Fix $B\in T(f)$. Let
$\cD_B(u,v)=\prod_j (v-f(u,\gamma_j(u)))$
where the product runs over all $j$ such that $\gamma_j$ leaves $T(f)$ at $B$. Set
$\cD_{\classB}(u,v)=\prod_{B'\in\classB}\cD_{B'}(u,v)$. 
Then $\cD_{\classB}(u,v)$ is a polynomial in $v$ with coefficients 
in $\bC\{u^{1/D}\}$. Furthermore we have:

\begin{Lemma}
$\cD_{\classB}(u,v)\in \bC\{u\}[v]$.
\end{Lemma}

\noindent \begin{proof}
 It is enough
to verify that for every complex number $v_0$ the index of
$\cD_{\classB}(u,v_0)\in \bC\{u^{1/D}\}$ is  1, which is equivalent, by Property \ref{index}, that the action
of  $\bU_D$ on this Puiseux series is trivial. 

\medskip
\noindent Take $\theta \in \bU_D$ and $B'\in \classB$.  We have 
$\theta * \cD_{B'}(u,v_0)=\prod_{j}
(v_0-f(u,\theta *\gamma_j(u))),
$
\noindent where $j$ runs over $\gamma_j$ leaving $T(f)$ at $B'$ and 
$
 \cD_{\theta *B'}(u,v_0)=\prod_{j}
(v_0-f(u,\gamma_j(u))),
$

\noindent where $j$ runs over $\gamma_j$ leaving $T(f)$ at $\theta*B'$.

\medskip
\noindent  Since $\gamma\in \Zer \frac{\partial f}
{\partial y}$ leaves $T(f)$ at $B'$ if and only if $\theta * \gamma$ leaves $T(f)$ at
$\theta * B'$, we get $\theta * \cD_{B'}(u,v_0)=\cD_{\theta *B'}(u,v_0)$. As a consequence
$\theta *\cD_{\classB}(u,v_0)=\theta * \prod_{B'\in \classB}\cD_{B'}(u,v_0)=\prod_{B' \in \classB}\cD_{\theta* B'}(u,v_0)=
\cD_{\classB}(u,v_0).$
\end{proof}

\medskip

\noindent  We conclude that $\prod_{\classB\in
E(f)}D_{\classB}(u,v)$ is an analytical factorization (not necessarily into irreducible factors) of the discriminant.

\medskip

\noindent By Proposition~\ref{invariant} every factor $\cD_{\classB}(u,v)$ has an elementary
Newton diagram
 of inclination $q(B)$. Observe that if
$\cD_{\classB}(u,v)$ is degenerate  then $\cD(u,v)$
is also degenerate. The aim of this section is to compute
the initial Newton polynomial of $\cD_{\classB}(u,v)$ . For this we need the next auxiliary results:

\begin{Lemma}
\label{arithmetical} Let $A,B$ be positive integers.  Then
$$\prod_{\theta^A=1}(z-\theta^Ba)=\left(z^{A/\gcd(A,B)}-a^{A/\gcd(A,B)}\right)^{\gcd(A,B)}.$$
\end{Lemma}

\noindent \begin{proof}Set $C=\gcd(A,B)$ and $A_ 1=A/C$,  $B_1=B/C$. Then
\begin{eqnarray*}
\prod_{\theta^A=1}\left(z-\theta^Ba\right)&=&\prod_{(\theta^C)^{A_1}=1}
\left(z-(\theta^C)^{B_1}a\right)=
\prod_{\small{\begin{array}{c}\omega^{A_1}=1\\ \theta^C =\omega \end{array}}}
\left(z-\omega^{B_1}a\right)\\ &=&
\prod_{\omega^{A_1}=1}\left(z-\omega^{B_1}a\right)^C=\left(z^{A_1}-a^{A_1}\right)^C,
\end{eqnarray*}

\noindent where the last equality holds
since the numbers $\omega^{B_1}a$ for $\omega^{A_1}= 1$ are all $A_1$-th
complex roots of $a^{A_1}$.
\end{proof}

\medskip

\begin{Lemma}\label{group}
Let $G$ be a finite group and $A$ be a finite set. Assume that $G$
acts on $A$ transitively, that is $A=Ga_0$ for some $a_0\in A$. Let
$P$ be a complex valued function on $A$. Set $G_0:=\{g\in
G\;:\;ga_0=a_0\}$. Then
\begin{enumerate}

\item $\sharp A \cdot \sharp G_0=\sharp G$.
\item $\prod_{g\in G}P(ga_0)=\prod_{a\in
A}\left(P(a)\right)^{\sharp G_0}$.
\end{enumerate}
\end{Lemma}

\noindent \begin{proof}  The first statement is the orbit-stabilizer theorem.

\noindent To prove the second statement consider the function $h:G\longrightarrow A$
given by $h(g)=ga_0$. Then $\prod_{g\in G}P(ga_0)=\prod_{a\in
A}\prod_{g\in h^{-1}(a)}P(h(g))=\prod_{a\in A}P(a)^{\sharp G_0}$.
The last equality holds
since the fibers of the function $h$ are the left-cosets of $G_0$ in $G$.
\end{proof}

\medskip

\noindent Now, our aim is to give a formula for $F_B(z)$ from Lemma \ref{L2}.

\medskip

\noindent Fix a pseudo-ball $B$ of $T(f)$. Let $f=f_1\cdots f_r$  be the decomposition of $f$ into
irreducible factors. Assume that $\Zer f_j\cap B\neq \emptyset$ for
$j\in \{1,\ldots,s\}$ and $\Zer f_j\cap B= \emptyset$ for $j\in
\{s+1,\ldots,r\}$. Note that $s\geq 1$ and perhaps $s=r$.
For every $j\in \{1,\ldots,s\}$ choose a  Newton--Puiseux root of $f_j(x,y)$  of the form
\begin{equation}
\label{truncation}
\lambda_B(x)+c_j x^{h(B)}+\cdots .
\end{equation}

\noindent Let $N$ be the index of $\lambda_B$ and write
$h(B)=\frac{m}{nN}$ with $m,n$ coprime.

\begin{Formula}
\label{Formula}
Keeping the above notations we have
\[F_B(z)=C \prod_{j=1}^s(z^n-c_j^n)^{\frac{\ord f_j(0,y)}{nN}}\]

\noindent where $C$ is a nonzero constant.
\end{Formula}

\noindent \begin{proof}
Fix $j\in \{1,\ldots,s\}$ and a Newton--Puiseux root $\alpha(x)$ of $f_j(x,y)$ of the form~(\ref{truncation}).
Since $f_j(x,y)$ is irreducible, the orbit $\bU_D\ast \alpha$ is the set $\Zer f_j$.
By  Lemma~\ref{group} the stabilizer $G_0$ of $\alpha(x)$ has $D/(\sharp \Zer f_j)=D/\ord f_j(0,y)$ elements.
Since every subgroup of a finite cyclic group is determined by the number of its elements, $G_0=\bU_{D/\ord f_j(0,y)}$.

\medskip

\noindent Let us observe that $\theta\ast\alpha$ belongs to $B$ if and only if $\theta\ast\lambda_B=\lambda_B$.
By a similar argument as before, the stabilizer $G_1$  of $\lambda_B$
is the subgroup~$\bU_{D/N}$ of~$\bU_D$.  Hence $\Zer f_j\cap B=G_1 \ast \alpha$.
By (ii) of Lemma~\ref{group} we get

\begin{equation}
\label{triangle}
\prod_{\theta \in G_1}(z-\lc_B(\theta \ast \alpha))=\prod_{\alpha_i \in
\Zer f_j \cap B }(z-\lc_B(\alpha_i))^{\frac{D}{\ord f_j(0,y)}}.
\end{equation}

\noindent On the other hand, following  Lemma \ref{arithmetical} we have
\begin{equation}
\label{square}
\prod_{\theta \in G_1}(z-\lc_B(\theta \ast \alpha))=\prod_{\theta^{D/N}=1}
(z-c_j\theta^{h(B)D})=(z^n-c_j^n)^{D/nN}.
\end{equation}

\noindent Comparing $(\ref{triangle})$ and $(\ref{square})$ we get
$\prod_{\alpha_i \in
\Zer f_j \cap B }(z-\lc_B(\alpha_i))=
(z^n-c_j^n)^{\ord f_j(0,y)/nN}.$

\noindent Finally $F_B(z)=
C\displaystyle \prod_{j=1}^s\prod_{\alpha_i\in \Zer f_j\cap B}(z-\lc_{B}(\alpha_i))=C\displaystyle
\prod_{j=1}^s (z^n-c_j^n)^{\ord f_j(0,y)/nN}$.
\end{proof}

\medskip

\noindent
From now on up to the end of this section we fix $B\in T(f)$ and 
put $q(B)=\frac{L}{M}$ with $L,M$ coprime.

\medskip

\noindent Let $\frac{d}{dz}F_B(z)=C'(z-z_1)\cdots (z-z_l)$.
Set $w_i=F(z_i)$ for $1\leq i\leq l$ and let $I:=\{i\in \{1,\ldots,l\}\;:\;w_i\neq 0\}$. Keeping this  notation we have:

\begin{Lemma}
\label{weighted}
The initial Newton  polynomial of $\cD_B(u,v)$ is
$$\inN \cD_B(u,v)=\prod_{i\in I}\bigl(v-w_iu^{q(B)}\bigr).$$

\end{Lemma}

\noindent \begin{proof}
By Lemma \ref{L2} the initial Newton polynomial of $\cD_B(u,v)$ is equal to $\prod_j\left(v-F_B(\lc_B \gamma_j)u^{q(B)}\right)$
where the product runs over $j$ 
such that $\gamma_j$ leaves $T(f)$ at $B$. It follows from Lemmas \ref{L1} and \ref{caract} that the above product equals 
$\prod_{i\in I}\left(v-w_iu^{q(B)}\right)$.
\end{proof}

\medskip

\noindent
\begin{Proposition}
\label{initial-disc}
Let $f(x,y)=0$ be a reduced complex plane curve.
Take a pseudo-ball $B$ of $T(f)$ such that $q(B)=\frac{L}{M}$ with $L,M$ coprime. Let $N$ be the index of $\lambda_B$. Then
\begin{equation}\label{eq:D}
\inN \cD_{\classB}(u,v)=\prod_{i\in I}\left(v^{M}-w_i^{M}
u^{L}\right)^{N/M}.
\end{equation}
\end{Proposition}

\noindent \begin{proof}
Recall that $\classB$ is the orbit of $B$ under the $\ast$ action of the group $\bU_D$.
Since $\theta \ast B=B$ if and only if $\theta \ast \lambda_B=\lambda_B$, the stabilizer
of $B$ is the subgroup $\bU_{D/N}$ (see the proof of Formula \ref{Formula}).

\medskip

\noindent We claim that under the assumptions
of Lemma \ref{weighted} one has $\inN \cD_{\theta\ast B}(u,v)=\prod_{i\in I} (v-w_i\theta^{q(B)D}u^{q(B)})$.
Indeed, by Proposition \ref{invariant} the critical values of $F_{\theta \ast B}$
are the critical values of $F_B$ times $\theta^{q(B)D}$, which proves the claim.

\medskip

\noindent By $(ii)$ of Lemma \ref{group} we have
\[
\prod_{\theta \in \bU_D}\inN \cD_{\theta \ast B}(u,v)=\prod_{B'\in \classB}\inN \cD_{B'}(u,v)^{D/N}=
\inN \cD_{\classB}(u,v)^{D/N}.
\]

\noindent On the other hand, by the claim  and Lemma \ref{arithmetical} we have

\[
\prod_{\theta \in \bU_D}\inN \cD_{\theta \ast B}(u,v)=\prod_{\theta^{D}=1}\prod_{i\in I}\left(v-w_i\theta^{q(B)D}u^{q(B)}\right)=\prod_{i\in I}\left(v^M-w_i^Mu^{L}\right)^{D/M}.
\]

\noindent Comparing the above equalities we get the proposition.
\end{proof}

\section{The irreducible case}
\label{irreducible}

\noindent We assume in this section that $f(x,y)\in \mathbf
C\{x,y\}$ is irreducible. Let $p:=\ord_yf(0,y)>1$ and $\Zer
f=\{\alpha_i(x)\}_{i=1}^{p}$.
The contacts $\{O(\alpha_i,\alpha_j)\}_{i\neq j}$, called {\em{the characteristic exponents of $f(x,y)$}}, form a finite set of rational numbers
$\left\{\frac{b_k}{p}
\right\}_{k=1}^h$, where $b_1<\ldots<b_h$. Set $b_0=p$. The sequence $(b_0,b_1,\ldots,b_h)$ is named {\em Puiseux characteristic}. Since $f(x,y)$ is irreducible, its Newton--Puiseux roots conjugate and all the pseudo-balls with the same height belong
to the same conjugate class in $E(f)$. Write $E(f)=\left\{\classB_1,\ldots,\classB_h\right\}$, where  $h(B_k)=\frac{b_k}{p}$ for 
$k\in\{1,\ldots,h\}$. By Lemma \ref{successor} $q(B_1)<q(B_2)<\cdots<q(B_h)$. The discriminant $\cD(u,v)$ is the product $\prod_{k=1}^h\cD_{\classB_k}(u,v)$.

\medskip

\noindent We now characterize the factors appearing in this product. Let $B\in T(f)$. By Formula \ref{Formula}, we have 
$F_B(z)=C (z^n-c^n)^{\frac{p}{nN}}$. This polynomial has only one nonzero critical value $w=F_B(0)$ of
multiplicity $n-1$. By Proposition \ref{initial-disc}, we have
$\inN \cD_{\classB}(u,v)=\left(v^{M}-w^{M}
u^{L}\right)^{(n-1)N/M}$, where $q(B)=\frac{L}{M}$,  $\gcd(L,M)=1$ and $N$ is the index of $\lambda_B$.
We stress that in the next corollary we only use the fact that $ \inN \cD_{\classB}(u,v)$ is a power of a quasi-homogeneous irreducible polynomial.

\medskip

\begin{Corollary} 
\label{no lattice points11} The power series $\cD_{\classB_i}(u,v)$  is non-degenerate  if and
only if there are no lattice points inside the only compact edge of its Newton diagram.
\end{Corollary}

\noindent Theorem \ref{teorema-ramas} is a consequence of Corollary \ref{no lattice points11} since the Newton diagram of $\cD(u,v)$ is the sum of the elementary Newton diagrams of 
$\cD_{\classB_i}(u,v)$.

\medskip

\noindent According to Merle \cite{Merle} and
Ephraim \cite{Ephraim}  the {\em semigroup}  $\Gamma$  (see for example \cite{Zariski} in the transversal 
case and \cite{GLP} in the general case) of 
$f(x,y)=0$  admits the minimal sequence
of generators $ \overline{b_0}:=\ord f(0,y),\overline{b_1}<\ldots<
\overline{b_h}$ and the Newton diagram of the discriminant $\cD(u,v)$ is

\begin{equation}
\label{dddd}
\sum_{k=1}^h \left\{\Teisssr{(n_k-1)\overline{b_k}}
{n_1\cdots n_{k-1}(n_k-1)}{18}{9}\right \},
\end{equation}

\noindent where $n_k:=\frac{\gcd(\overline{b_0},\overline{b_1},\ldots,
\overline{b_{k-1}})}{\gcd(\overline{b_0},\overline{b_1},
\ldots,\overline{b_k})}=\frac{\gcd(b_0,\ldots,b_{k-1})}{\gcd(b_0,\ldots,b_{k})}$ and by convention $n_0=1$.
The inclinations of the edges of the Newton diagram (\ref{dddd}) are $q(B_1),\ldots, q(B_h)$. They are called {\em polar invariants} of the pair $(x,f)$.

\medskip

\noindent Since the Newton diagram of a product is the sum of the Newton diagrams of its factors and
the sequence $(q(B_k))$ is increasing, the Newton diagram of $\cD_{\classB_k}(u,v)$ is the $k$-th term of (\ref{dddd}).

\begin{Corollary}
\label{coro-irre}
The power series $\cD_{\classB _k}(u,v)$ is
non-degenerate if and only if $(n_k-1)\gcd(\overline{b_k},n_1\cdots n_{k-1})=1$.
\end{Corollary}

\noindent \begin{proof}
Since the Newton diagram of 
 $\cD_{\classB _k}(u,v)$ is  $\left\{\Teisssr{(n_k-1)\overline{b_k}}
{n_1\cdots n_{k-1}(n_k-1)}{18}{9}\right \}$ the statement follows from Corollary \ref{no lattice points11}.
\end{proof}

\medskip

\begin{Remark}
Note that if for $k>1$ the polar invariant $q(B_k)$ is an integer then  $\left\{\Teisssr{(n_k-1)\overline{b_k}}
{n_1\cdots n_{k-1}(n_k-1)}{18}{9}\right \}$ has lattice points inside its compact edge and $\cD(u,v)$ is degenerate.
\end{Remark}

\noindent Observe that a necessary condition for $\cD(u,v)$ to be non-degenerate is $n_1=n_2=\ldots=n_h=2$, where $h$ is the number of characteristic exponents of $f=0$.

\medskip
\begin{Corollary}
\label{caract-branches}
Let $f(x,y)=0$ be a branch with  $h$ characteristic exponents.
We have
\begin{enumerate}
\item If $h=1$ then the discriminant $\cD(u,v)$ is non-degenerate if and only if $\ord f(0,y)=2$.
\item If $h=2$ then  the discriminant $\cD(u,v)$ is non-degenerate if and only if $\ord f(0,y)=4$.
\item If $h>2$ then $\cD(u,v)$ is degenerate.
\end{enumerate}
\end{Corollary}

\section{The general case}
\label{general-section}
\noindent In this section we specify the polynomial factorization of $\inN \cD_{\classB}(u,v)$. 
We start with four technical lemmas. 
Their sole purpose is to show that the factors of~(\ref{equatt1}) and~(\ref{equat2}) 
in Proposition~\ref{polynomial-factorization}  are polynomials.

\begin{Lemma}
\label {LEMMAN}
Let $0\leq a\leq b$ and let $f:[a,b]\to {\bf R}$ be a continuous function 
such that $f(x)\geq 0$ for $a\leq x\leq b$. 
Let $c$ be a positive integer. 
Then $\max_{x\in[a,b]} {{(m-x)^c}\over{m^c}}f(x) \to \max_{x\in[a,b]} f(x)$ as $m\to \infty$.
\end{Lemma}

\noindent
\begin{proof}
Let $x_0$ be the point of the interval $[a,b]$ such that $f(x_0)=\max_{x\in[a,b]} f(x)$.  
We have  ${{(m-x_0)^c}\over{m^c}}f(x_0)  \leq 
\max_{x\in[a,b]} {{(m-x)^c}\over{m^c}}f(x)\leq 
\max_{x\in[a,b]} f(x)$ for large $m$. Passing to the limits we get the lemma. 
\end{proof}

\begin{Lemma}
\label{Morse1}
Let $a_1$, \dots , $a_n$ be positive integers. 
Then there exist  pairwise different nonzero complex numbers $d_1$, \dots $d_n$ such that the polynomial 
$$ H(t)=\prod_{j=1}^n(t-d_j)^{a_j}
$$
has $n-1$ pairwise different nonzero critical values, and all of them differ from $H(0)$.
\end{Lemma}

\noindent
\begin{proof}
It suffices to construct step-by-step a sequence  $0<d_1<d_2<\dots <d_n$ 
such that the polynomials $W_k(t)=\prod_{j=1}^k(t-d_j)^{2a_j}$ for $k\in \{1,\dots, n\}$ 
satisfy conditions $W_k(0)<\max_{t\in[d_1,d_2]}W_k(t)<\dots <\max_{t\in[d_{k-1},d_k]}W_k(t)$. 

\noindent Assume that the numbers $0<d_1<\dots<d_k$ and the polynomial $W_k(t)$ are already 
constructed. Applying Lemma \ref{LEMMAN} 
to every interval $[d_{j-1},d_j]$ and to the interval $[0,0]$ we conclude that for sufficiently large  $m=:d_{k+1}$ the maximal values 
of the polynomial 
$$
{{1}\over{m^{2a_{k+1}}}}W_{k+1}(t)={{(m-t)^{2a_{k+1}}}\over{m^{2a_{k+1}}}}W_k(t)
$$
in the intervals $[0,0]$, $[d_1,d_2]$, \dots, $[d_{k-1},d_k]$ form an increasing sequence and  are bigger than $W_{k+1}(0)/m^{2a_{k+1}}$.

\noindent To assure that $\max_{t\in[0,d_k]}W_{k+1}(t)<\max_{t\in[d_k,d_{k+1}]}W_{k+1}(t)$
it is enough to observe that in the sequence of inequalities
\begin{eqnarray*}
\max_{t\in[0,d_k]}W_{k+1}(t)& \leq&
m^{2a_{k+1}}\max_{t\in[0,d_k]}W_{k}(t) <
\left({{m-d_k}\over{2}}\right)^{\deg W_{k+1}(t)} \\
& \leq &
W_{k+1}\left({{m+d_k}\over{2}}\right) \leq 
\max_{t\in[d_k,d_{k+1}]}W_{k+1}(t),
\end{eqnarray*}

\noindent the second inequality holds for all $m$ big enough. Finally taking $H(t):=\prod_{j=1}^n(t-d_j)^{a_j}$ we see that the nonzero critical values of $W_n(t)$ are the squares of the nonzero critical values of $H(t)$ and we prove the lemma.
\end{proof}

\begin{Corollary}
\label{Morse2}
Let $H(t)$ be a complex polynomial of the form 
\begin{equation}\label{Eq:H1}        
 H(t)=t^{a_0}\prod_{j=1}^{n}(t-d_j)^{a_j} ,
\end{equation}
where  $a_j$ are positive integers for $j\in \{0,1,\ldots ,n\}$.
Then for some $d_1$,\dots, $d_n$ the polynomial $H(t)$ has $n$  pairwise different nonzero critical values.
\end{Corollary} 

\noindent
\begin{proof}
By Lemma \ref{Morse1} we can choose a sequence $e_0,e_1,\ldots,e_n$ such that the polynomial $H_1(t)=\prod_{j=0}^n(t-e_j)^{a_j}$ has $n$ pairwise different nonzero critical values. We finish by putting $H(t)=H_1(t+e_0)$ and $d_j=e_j-e_0$ for $j\in\{1,\ldots,n\}$.
\end{proof}

\medskip
\noindent
In the next lemma we change the notation slightly. 
Notice that the polynomial~$F_B(z)$ and the exponent~$q(B)$ in Lemma \ref{L2} 
depend not only on $B$ but also on the power series $f(x,y)$. We write 
$F_{B,f}(z)$ for the polynomial and $q(B,f)$ for the exponent to stress this dependence. 

\begin{Lemma}
\label{change_of_series}
Let $f(x,y)$ be a reduced power series such that $f(0,y)\neq0$. 
Fix $B\in T(f)$. Let $N$ be the index of $\lambda_B$ and write $h(B)=\frac{m}{nN}$ with $m$, $n$ coprime. 
Assume that 
$ F_{B,f}(z)=Cz^{a_0}\prod_{j=1}^s(z^n-d_j)^{a_j}$, 
where $d_j$ are pairwise different nonzero complex numbers, $a_0$ is a nonnegative integer and
$a_j$ are positive integers for $j\in\{1,\dots,s\}$.

\noindent Then for every sequence of pairwise different nonzero complex numbers 
$\tilde d_1$,\dots, $\tilde d_s$  there exists a reduced power series $\tilde f(x,y)$ such that 
$B\in T(\tilde f)$, $q(B,\tilde f)=q(B,f)$ and 
$ F_{B,\tilde f}(z)=Cz^{a_0}\prod_{j=1}^s(z^n-\tilde d_j)^{a_j}$.
\end{Lemma} 

\noindent
\begin{proof}
Let $f=f_1\cdots f_r$  be the decomposition of $f(x,y)$ into irreducible factors. 
Without loss of generality we may assume that $\Zer f_i\cap B\neq \emptyset$ for
$i\in \{1,\ldots,k\}$ and $\Zer f_i\cap B= \emptyset$ for $i\in\{k+1,\ldots,r\}$. 
For every $i\in \{1,\ldots,k\}$ choose a  Newton--Puiseux root of $f_i$  of the form
$\alpha_i(x)=\lambda_B(x)+c_i x^{h(B)}+\cdots\;$.
Let $\cC=\{c_i^n\;:\;i\in\{1,\dots,k\}\}$. Then it follows from Formula~\ref{Formula} that 
$\cC\setminus\{0\}=\{d_1,\dots,d_s\}$, 
$a_0=\frac{1}{N}\sum_{i:c_i=0}\ord f_i(0,y)$ and 
$a_j=\frac{1}{nN}\sum_{i:c_i^n=d_j}\ord f_i(0,y)$ for $j=1,\dots,s$.

\noindent For every $i\in \{1,\ldots,k\}$ take the fractional power series
$$ \tilde\alpha_i(x)=\alpha_i(x)+(\tilde c_i-c_i)x^{h(B)}=\lambda_B(x)+\tilde c_i x^{h(B)}+\cdots
$$
where $\tilde c_i=0$ if $c_i=0$ and $\tilde c_i^n=\tilde d_j$ if  $c_i^n=d_j$.
Set $\tilde f=a\tilde f_1\cdots \tilde f_k f_{k+1}\cdots f_r$,  where 
$\tilde f_i(x,y)$ are irreducible power series such that $\tilde\alpha_i\in \Zer\tilde f_i$  for $i\in\{1,\ldots,k\}$ and 
$a$ is a constant which will be specified later.  
Clearly $B$ is an element of $T(\tilde f)$. 

\noindent Now let us compute $F_{B,\tilde f}$.  One has $\ord f_i(0,y)=\ord \tilde f_i(0,y)$ for $i=1,\ldots,k$ since 
$\alpha_i(x)$ and $\tilde \alpha_i(x)$ have the same index. By the first part 
of the proof it is clear that $F_{B,\tilde f}(z)=\tilde Cz^{a_0}\prod_{j=1}^s(z^n-\tilde d_j)^{a_j}$. 
By a suitable choice of the complex number $a$ we get~$\tilde C=C$. 

\medskip

\noindent 
It remains to prove that $q(B,f)=q(B,\tilde f)$. 
Let $\gamma(x)=\lambda_B(x)+cx^{h(B)}$ where~$c$ is a generic constant. 
Then $q(B, f)=\ord f(x,\gamma(x))=\sum_{i=1}^r \ord f_i(x,\gamma(x))$ and an 
analogous formula holds for $q(B, \tilde f)$. 

\medskip
\noindent Fix $i\in\{1,\dots,k\}$.
For generic $c$ we have $\cont(\gamma,\Zer f_i)=\cont(\gamma,\Zer \tilde f_i)=h(B)$. Since the Puiseux 
characteristics of both irreducible power series are the same, we get $\ord f_i(x,\gamma(x))=\ord\tilde f_i(x,\gamma(x))$
(see for example \cite{Merle}, Proposition 2.4 for the transverse case and \cite{GwP}, Proposition 3.3 for the general case).
\end{proof}

\medskip
\noindent
\textbf{Remark.}
One can show that the power series $\tilde f(x,y)$ constructed in the proof of Lemma \ref{change_of_series} has the same equisingularity 
type as~$f(x,y)$.

\medskip
\noindent We introduce a new  polynomial $H_B(t)$ 
associated with $B\in T(f)$ whose critical values provide  a polynomial factorization of $\inN D_{\classB}(u,v)$.

\medskip
\begin{Lemma}
Fix $B\in T(f)$. Let  $N$ be the index of $\lambda_{B}$. Write
$h(B)=\frac{m}{nN}$  and $q(B)=\frac{L}{M}$ where $\gcd(m,n)=\gcd(L,M)=1$.  
Then there exists a unique polynomial $H_B(t)$ such that 
$H_B(z^{n})=F_B(z)^M.$
\end{Lemma}
\noindent \begin{proof}
\noindent Assume as earlier that all Newton--Puiseux roots of $f(x,y)$
and $\frac{\partial f}{\partial y}(x,y)$
belong to $\bC\{x^{1/D}\}$ for some positive integer $D$. We use the properties of the conjugate action introduced in Section \ref{section-factorization}. One easily checks that $\theta*B=B$ for $\theta \in \bU_{D/N}$ (see the proof of Proposition \ref{initial-disc}). Set $D=D_0nN$ and take $\theta \in \bU_{D/N}$ such that $\omega:=\theta^{D_0}$ is an $n$-th primitive root of unity.  By
Proposition \ref{invariant} we get $\theta^{q(B)D}F_B(z)=F_{B}(\theta^{h(B)D}z)$. Hence
$F_B(z)^M=F_B(\omega^mz)^M$. Comparing the terms  of both sides we see that all monomials appearing in the polynomial $F_B(z)^M$ are powers of $z^n$.
\end{proof}

\medskip

\begin{Proposition}\label{polynomial-factorization}
Let $f(x,y)=0$ be a reduced curve. Fix $B\in T(f)$.  
Let  $N$ be the index of $\lambda_{B}$. Write
$h(B)=\frac{m}{nN}$ and 
$q(B)=\frac{L}{M}$  where $\gcd(m,n)=\gcd(L,M)=1$.
Let $H'_B(t)= C (t-t_1)\cdots (t-t_r)$. Set ${\textsf w}_0=H_B(0)$, ${\textsf w}_j=H_B(t_j)$ and
$J=\{j\in\{1,\ldots,r\}\;:\; {\textsf w}_j\neq 0\}$. Then

\begin{eqnarray}
 \inN \cD_{\classB}(u,v)=(v^{M}-{\textsf w}_0u^{L})^{(n-1)N/M}\prod_{j\in J}(v^{M}-{\textsf w}_ju^{L})^{nN/M}
& \mbox{if ${\textsf w}_0\neq0$,} &
  \label{equatt1}\\
\inN \cD_{\classB}(u,v)=\prod_{j\in J}(v^{M}-{\textsf w}_ju^{L})^{nN/M}
& \mbox{if ${\textsf w}_0=0$.} &
\label{equat2}
\end{eqnarray}
Moreover (\ref{equatt1}) and (\ref{equat2}) give a polynomial factorization of $\inN \cD_{\classB}(u,v)$.
\end{Proposition}

\noindent
\begin{proof}

\noindent
The above formulas follow from Proposition~\ref{initial-disc} and the equality 
$MF_B(z)^{M-1}F_B'(z)$\\$=nz^{n-1}H_B'(z^{n})$ 
 which allows to express critical values of $F_B$ in terms of critical 
values of $H_B$. 

\noindent Using Lemma~\ref{change_of_series} we can replace $f(x,y)$ by such a power series $\tilde f(x,y)$ that 
conclusions of  Lemma~\ref{Morse1} or Corollary~\ref{Morse2}, for $H(t)=H_B(t)$,  are satisfied. 
Then $\{{\textsf w}_j\}_{j\in J\cup \{0\}}$ is a sequence of pairwise different complex numbers. 
The polynomials $v^{M}-{\textsf w}_ju^{L}$ are irreducible and pairwise coprime. Hence the exponents $(n-1)N/M$, $nN/M$ in~(\ref{equatt1})  or $nN/M$ in~(\ref{equat2}) are integers.
\end{proof}

\begin{Theorem}\label{Th:general}
Let $f(x,y)=0$ be a reduced curve and let $B\in T(f)$.  
Let  $N$ be the index of $\lambda_{B}$. Write
$h(B)=\frac{m}{nN}$ and 
$q(B)=\frac{L}{M}$  where $\gcd(m,n)=\gcd(L,M)=1$.

\begin{enumerate}
\item If $H_B(t)$ has only one root (possibly multiple), then $\cD_{\classB}(u,v)$ is non-degenerate if and only if $(n-1)N=M$.
    
\item Otherwise  $\cD_{\classB}(u,v)$ is non-degenerate if and only if $nN=M$ and all
nonzero critical values of $H_B(t)$ are simple.
\end{enumerate}
\end{Theorem}

\noindent 
\begin{proof}

\noindent Assume that $H_B(t)$ has only one root. 
By Proposition~\ref{polynomial-factorization}
$\inN \cD_{\classB}(u,v)=(v^{M}-{\textsf w}_0u^{L})^{(n-1)N/M}$. 
This polynomial is non-degenerate if and only if $(n-1)N=M$.

\medskip
\noindent 
Suppose that $H_B(t)$ has at least two different roots. Assume that ${\textsf w}_0=0$. Then (\ref{equat2}) is a reduced polynomial if and only if $nN=M$ and all
nonzero critical values of $H_B(t)$ are simple. Assume now that ${\textsf w}_0\neq 0$. Then the polynomial (\ref{equatt1}) is reduced  if and only if
$nN/M=1$ and $({\textsf w}_j)_{j\in J}$ is a sequence 
of pairwise diffe\-rent complex numbers. 
In this case the only difficulty  arrives from the term $(v^{M}-{\textsf w}_0u^{L})^{(n-1)N/M}$ but the exponents  $nN/M$ and $(n-1)N/M$ are integers, so
the condition $nN/M=1$ 
 implies $(n-1)N/M=0$.
\end{proof}

\medskip

\noindent We finish this section with another example  of a multibranched curve $f=0$ such that the discriminant of the morphism $(x,f)$ is non-degenerate. For the construction we use the {\em Eggers tree} whose construction we now recall. 
We assume that $x=0$ and $f=0$ are transverse.
\noindent Recall that $E(f)$ is the set of all conjugate classes of $B$ for $B\in T(f)$.  An element of $E(f)$ is uniquely determined 
by its height $h(\classB):=h(B)$ and the set of irreducible factors $f_i$ of $f$ such that $\Zer f_i\cap B\neq\emptyset$
(see \cite{Kuo-Parusinski2}, Section 6). The tree structure on $T(f)$ induces a tree structure on $E(f)\cup\{\,f_0,\dots,f_k\,\}$.
This newly constructed tree is called the Eggers tree of $f$ (\cite{Eggers}, see also \cite{GB}). In Eggers' terminology the vertices from $E(f)$ are called {\em black points} and those from $\{\,f_0,\dots,f_k\,\}$ are called {\em white points}. The Eggers tree is an oriented tree where the root is the black point of the minimal height and the leaves are the white points. The {\em outdegree} of a vertex $Q$  is the number of edges joining $Q$ with its successors.

\begin{Remark}
\noindent The first part in Theorem \ref{Th:general} corresponds to simple points (i.e. vertices of outdegree 1)  in the Eggers tree. The second part corresponds to bifurcation points (vertices of outdegree greater than 1)  in the Eggers tree. The number of  irreducible factors of $H_B(t)$ is equal to the outdegree of the vertex $\classB$.
\end{Remark}

\begin{Example}{\rm 
Set $n_0=1$ and let $n_1, \dots, n_k$  be pairwise coprime integers bigger than 1. 
We construct a singular power series $f=f_0f_1\cdots f_k$, where $f_i$ are irreducible power series, 
$\ord f_i(0,y)=n_0\cdots n_i$ for $i\in\{0,\dots, k\}$, and such that the discriminant of the morphism $(x,f)$ is non-degenerate.

\medskip

\noindent Let $h_i=1+\frac{1}{n_1}+\cdots+\frac{1}{n_i}$ for $i\in\{1,\dots,k\}$. 
We claim that $h_i$ can be written as $\frac{b_i}{n_1\cdots n_i}$, 
with  $b_i$ and $n_1\cdots n_i$  coprime.
The proof runs by induction on $i$. For $i=1$ we have $h_1=\frac{n_1+1}{n_1}$.  Assume that $\gcd(b_i,n_1\cdots n_i)=1$. 
By equality $\frac{b_{i+1}}{n_1\cdots n_{i+1}}=\frac{b_i}{n_1\cdots n_i}+\frac{1}{n_{i+1}}$ we get 
$b_{i+1}=b_in_{i+1}+n_1\cdots n_i$. 
Thus $\gcd(b_{i+1},n_{i+1})=\gcd(n_1\cdots n_i,n_{i+1})=1$,
$\gcd(b_{i+1},n_1\cdots n_i)=\gcd(b_in_{i+1},n_1\cdots n_i)=1$ and consequently we get
$\gcd(b_{i+1},n_1\cdots n_{i+1})=1$.

\medskip

\noindent Let 
\begin{eqnarray*}
\alpha_0(x)&=& 0, \\
\alpha_1(x)&=&x^{h_1}, \\
\alpha_2(x)&=&x^{h_1}+x^{h_2}, \\
&\vdots& \\
\alpha_k(x)&=&x^{h_1}+x^{h_2}+\cdots+x^{h_k}. 
\end{eqnarray*}
\noindent We consider $f=f_0f_1\cdots f_k$ where $f_i$ are irreducible power series such that $\alpha_i\in\Zer f_i$. 
By Property \ref{index} the order of $f_i(0,y)$ is $n_0\cdots n_i$ for $i\in \{0,\dots, k\}$.
Let $B_i=B(\alpha_{i-1},h_i)$ for $i\in \{1,\dots,k\}$. 
Then $E(f)=\{\,\classB_1,\dots,\classB_k\,\}$.  The Eggers tree of $f$ is drawn below. 

\[
\includegraphics[width=4cm]{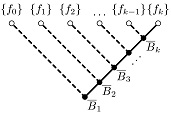}
\]

\noindent 
Since $\lambda_{B_i}(x)=\alpha_{i-1}(x)$  we have,
with the notations of  Formula \ref{Formula},
$N=n_0\cdots n_{i-1}$ and $n=n_i$. 
Hence 
\begin{equation}\label{Eq:F}
F_{B_i}(z)=
C(z^{n}-0)^{\frac{\ord f_{i-1}(0,y)}{nN}}
\prod_{j=i-1}^k(z^n-1)^{\frac{\ord f_j(0,y)}{nN}}=
Cz(z^{n_i}-1)^{A_i}, 
\end{equation}
where $A_i$ is a positive integer.

\medskip

\noindent Now we show that $q(B_i)$ could be written as $\frac{L_i}{M_i}=\frac{L_i}{nN}$ with $L_i$ and $nN$ coprime. 
Since $h(B_i)=\frac{b_i}{nN}$ with $b_i$ and $nN$  coprime,
it is enough to prove by induction on $i$ that for $i\in \{1,\dots,k\}$ the difference 
$q(B_i)-h(B_i)$ is an integer. By Lemma \ref{successor}  and (\ref{Eq:F}) we get
 $q(B_1)=\sharp(\Zer f)h(B_1)=\deg F_{B_1}(z)h(B_1)=(1+n_1A_1)h(B_1)$. Hence 
 $q(B_1)-h(B_1)=b_1A_1$. Now, again by ~(\ref{Eq:F}) and Lemma \ref{successor} we get
$q(B_{i+1})-q(B_i)=(1+n_{i+1}A_{i+1})\frac{1}{n_{i+1}}=\frac{1}{n_{i+1}}+A_{i+1}$.
Thus by the inductive hypothesis 
$q(B_{i+1})-h(B_{i+1})=q(B_i)+\frac{1}{n_{i+1}}+A_{i+1}-h(B_{i+1})=q(B_i)-h(B_i)+A_{i+1}$
is an integer. 

\medskip

\noindent The only roots of  $H_{B_i}(t)$ are $0$ and $1$. Therefore this polynomial has a unique nonzero critical value ${\textsf w}_i$. By equality $nN=M_i$ and Proposition~\ref{polynomial-factorization}  we get $\inN\,\cD_{\classB_i}(u,v)=v^{L_i}-{\textsf w}_iu^{M_i}$.

\medskip

\noindent The polynomials $\inN \cD_{\classB_i}(u,v)$, for $1\leq i\leq k$, are irreducible and pairwise coprime. Hence 
the discriminant $D(u,v)=\cD_{\classB_1}(u,v)\cdots \cD_{\classB_k}(u,v)$  of the morphism $(x,f)$ is non-degenerate.}
\end{Example}

\section {Stability of the discriminant's initial Newton polynomial}

\noindent To simplify subsequent statements we say that the power series $H_1(u,v)$, $H_2(u,v)$ are 
\emph{equal up to rescaling variables} if there exist nonzero constants $A$, $B$, $C$ such that 
$H_1(u,v)=CH_2(Au,Bv)$.  The Kouchnirenko non-degeneracy of a power series in two variables 
does not depend on rescaling variables. 

\begin{Lemma} 
\label{lc1}
Let $\cD(u,v)$ be the discriminant of the morphism $(f,g)$. Then for any  nonzero constants $A$, $B$  the discriminant  curve
of the morphism $(Af,Bg)$ admits the equation 
$\cD(u/A,v/B)=0$.
\end{Lemma}

\noindent \begin{proof} Let $u=Au'$, $v=Bv'$. 
As $(u,v)=(Af(x,y),Bg(x,y))$ then $(u',v')=(f(x,y),g(x,y))$. Hence, the discriminant curve of the morphism 
$(Af,Bg)$ admits the equation $\cD(u',v')=0$ which gives the lemma.  
\end{proof}

\begin{Theorem} 
\label{tc1}
Let $f=0$ be a reduced singular curve and let $\ell=0$ be  a smooth curve which is not a branch of $f=0$. 
Then for every invertible power series $u_1(x,y)\in \bC\{x,y\}$
the initial Newton  polynomials of the discriminants of  $(\ell,f)$ 
and $(\ell,u_1f)$ are equal up to rescaling variables. 
\end{Theorem}

\noindent \begin{proof}
An analytic change of coordinates does not affect the equation of the discriminant. Hence, we may assume 
that $\ell(x,y)=x$. By Lemma \ref{lc1} we may also assume that $u_1(0,0)=1$. 
Since $f$ and $u_1f$ have the same New\-ton--Puiseux roots, their tree models  coincide. Let $B\in T(f)$. Applying Lemma \ref{L2} to $f$ and $u_1f$ 
 we show that  $F_{B,f}(z)=F_{B,u_1f}(z)$ and $q(B,f)=q(B,u_1f)$. 
By Lemma~\ref{weighted} the initial Newton polynomial of the discriminant depends only on 
$F_{B}(z)$ and $q(B)$ for pseudo-balls $B$ from the tree model. This proves Theorem~\ref{tc1}
\end{proof}
\medskip

\noindent 
In what follows we need a few auxiliary results about fractional power series. 

\medskip

\noindent Consider the fractional power series $\phi(x)=x+{}\cdots=x(1+\phi_1(x))$.
We define the formal root $\phi(x)^{1/n}:=x^{1/n}\sqrt[n]{1+\phi_1(x)}$, where
$\sqrt[n]{1+z}:=1+\frac{1}{n}z+{}\cdots$ is an analytic branch of the $n$-th
complex root of $1+z$. Then, having a power series $\psi(x)=\overline\psi(x^{1/n})$,
where $\overline\psi(t)$ is a convergent power series, we define the formal
substitution $\psi(\phi(x))$ as the fractional power series
$\overline\psi\bigl(\phi(x)^{1/n}\bigr)$.
\begin{Lemma}
\label{lc2}
Let 
\begin{eqnarray*}
\alpha_i(x) &=& x+\sum_{k=n+1}^{N-1}a_kx^{k/n}+c_ix^{N/n}+{}\cdots \\
\beta_i(y) &=& y+\sum_{k=n+1}^{N-1}b_ky^{k/n}+d_iy^{N/n}+{}\cdots
\end{eqnarray*}
for $i=1,2$. 
If $\beta_1(\alpha_1(x))=\beta_2(\alpha_2(x))$ then $c_1-c_2=d_2-d_1$.
\end{Lemma}

\noindent \begin{proof}
Write $\lambda(y)=\sum_{k=n+1}^{N-1}b_ky^{k/n}$. 
Then 
\begin{eqnarray*}
0&=&\beta_1(\alpha_1(x))-\beta_2(\alpha_2(x))\\
& = &  [\alpha_1(x)-\alpha_2(x)] + 
[\lambda(\alpha_1(x))-\lambda(\alpha_2(x))] + 
[d_1(\alpha_1(x))^{N/n}-d_2(\alpha_2(x))^{N/n}]+{}\cdots \\ 
&= &  [(c_1-c_2)x^{N/n}+{}\cdots] + 
[(d_1-d_2)x^{N/n}+{}\cdots] +
[\lambda(\alpha_1(x))-\lambda(\alpha_2(x))] \\
&= &  [(c_1-c_2+d_1-d_2)x^{N/n}+{}\cdots]+
[\lambda(\alpha_1(x))-\lambda(\alpha_2(x))].
\end{eqnarray*}
To finish the proof it suffices to show that the fractional power series 
$\lambda(\alpha_1(x))-\lambda(\alpha_2(x))$ does not contain the term of 
order~$N/n$. This task reduces to 

\medskip
\noindent \textbf{Claim.} For every $k>n$ the order of  
$(\alpha_1(x))^{k/n}-(\alpha_2(x))^{k/n}$ is
 bigger than $N/n$.

\medskip
\noindent \textbf{Proof of the Claim.}  
For every convergent power series 
$g(z)\in\bC\{z\}$ there exists 
$G(z,w)\in\bC\{z,w\}$ such that $g(z)-g(w)=(z-w)G(z,w)$.

\medskip

\noindent Let $\alpha_i(x)=x(1+\tilde\alpha_i(x))$ for $i=1,2$. Using the above equality for $g(z)=\sqrt[n]{1+z}$ we get
$(\alpha_1(x))^{k/n}-(\alpha_2(x))^{k/n} = 
x^{k/n}\bigl(\bigl(\sqrt[n]{1+\tilde\alpha_1(x)}\bigr)^k- 
             \bigl(\sqrt[n]{1+\tilde\alpha_2(x)}\bigr)^k\bigr) = 
x^{k/n}(\tilde\alpha_1(x)-\tilde\alpha_2(x))G(\tilde\alpha_1(x),\tilde\alpha_2(x)) = 
x^{(k-n)/n}(\alpha_1(x)-\alpha_2(x))G(\tilde\alpha_1(x),\tilde\alpha_2(x))$
which proves the Claim. 
\end{proof}

\begin{Lemma} 
\label{lc3}
Let $f(x,y)=(y-x)^n+\cdots$ be an irreducible complex power series. 
Then for every Newton--Puiseux root $y=\alpha(x)$ of $f(x,y)$ 
there exists a Newton--Puiseux root $x=\beta(y)$ of $f(x,y)$ 
such that $\beta(\alpha(x))=x$. 
\end{Lemma}

\noindent \begin{proof}
Fix a Newton--Puiseux root   $y=\alpha(x)$ of $f(x,y)$.
Let $\beta_1(y),\dots,\beta_n(y)$ be  the solutions of $f(x,y)=0$
in $\bC\{y\}^*$.  Then there exists
a unit $v(x,y)\in\bC\{x,y\}$ such that
$f(x,y)=v(x,y)\prod_{j=1}^n(x-\beta_j(y))$.
By Property \ref{index} the index of every $\beta_j(y)$ is $n$ and we can 
write $\beta_j(y)=\bar\beta_j(y^{1/n})$.
Substituting $y:=s^n$ we get
$f(x,s^n)=v(x,s^n)\prod_{j=1}^n(x-\bar\beta_j(s))$.
By putting $s:=\alpha(x)^{1/n}$ we obtain
$$
0=f(x,\alpha(x))=v(x,\alpha(x))\prod_{j=1}^n(x-\bar\beta_j(\alpha(x)^{1/n}))
$$
and the lemma follows.
\end{proof}

\medskip
\begin{Remark}
By Lemma~\ref{lc3} for every fractional power series $y=\alpha(x)=x+\cdots$ there exists a fractional power series 
$x=\beta(y)$ such that  $\beta(\alpha(x))=x$. 
We call  $x=\beta(y)$ a {\em formal inverse} of $y=\alpha(x)$. 
By Lemma~\ref{lc2} a formal inverse is unique. One can also show that if 
$x=\beta(y)$ is the formal inverse of $y=\alpha(x)$ then 
$y=\alpha(x)$ is the formal inverse of $x=\beta(y)$. 
\end{Remark}

\begin{Theorem}
\label{tc3}
Let $f=0$ be a unitangent reduced singular curve and let $\ell_1=0$, $\ell_2=0$ be smooth curves 
transverse to $f=0$. Then the initial Newton polynomials of the dis\-criminants of morphisms $(\ell_1,f)$, 
$(\ell_2,f)$ are equal up to rescaling variables. 
\end{Theorem}

\noindent \begin{proof}
Assume that the curves $\ell_1=0$, $\ell_2=0$ are transverse. 
Then there exists a system of local analytic coordinates $(\tilde x,\tilde y)$ 
such that  $\ell_1=\tilde x$ and $\ell_2=\tilde y$.  
By assumption the curve $f=0$ has only one tangent $\tilde y=c\tilde x$, where $c\neq0$. 
In the new coordinates $(x,y)=(c\tilde x,\tilde y)$ this tangent becomes $y=x$.  

\medskip

\noindent  Let $g(x,y)$ be the Weierstrass polynomial of $f(x,y)$ 
and $g'(x,y)$ be the Weierstrass polynomial of $f(-y,x)$. 
Then by Lemma~\ref{lc1} and Theorem~\ref{tc1}  the initial Newton polynomials of the discriminants of the morphisms 
$(\ell_1,f)$ and $(x,g)$ are equal up to rescaling variables.
The same applies to the morphisms $(\ell_2,f)$ and $(x,g')$. 

\medskip

\noindent Write $\Zer g=\{\alpha_1(x),\ldots,\alpha_p(x) \}$.
Let $\beta_i(y)$ be  the formal inverse of $\alpha_i(x)$ for $i=1,\dots,p$.  
It follows from Lemma~\ref{lc3} that $\alpha_i'(x)=-\beta_i(x)$ are Newton-Puiseux roots of 
$g'(x,y)$ for $i=1,\dots,p$. 
By Lemma~\ref{lc2}  $\ini (\alpha_i(x)-\alpha_j(x))=\ini (\alpha'_i(x)-\alpha'_j(x))$ for $1\leq i<j\leq p$.
We get $\Zer g'=\{\alpha_1'(x),\ldots,\alpha_p'(x) \}$.

\medskip

\noindent
The mapping
$B(\alpha_i,O(\alpha_i,\alpha_j))\mapsto B(\alpha'_i,O(\alpha'_i,\alpha'_j))$ 
gives a one-to-one correspondence between pseudo-balls of the tree models $T(g)$ and $T(g')$. 
Moreover, for every $B\in T(g)$ and the corresponding $B'\in T(g')$ there exists a constant $a_B$ such 
that $\lc_{B'}(\alpha'_i)=\lc_{B}(\alpha_i)+a_B$ for $\alpha_i\in B$, $\alpha'_i\in B'$.  

\medskip

\noindent By Remark \ref{xxx}, the leading coefficients of $F_{B,g}(z)=C\prod_{i:\alpha_i \in B}(z-\lc_B(\alpha_i))$ and $F_{B',g'}(z)=C'\prod_{i:\alpha'_i \in B'}(z-\lc_{B'}(\alpha'_i))$ are given respectively by 
$$Cx^{q(B,g)}=\prod_{i:\alpha_i\not\in B}\ini (\alpha_j(x)-\alpha_i(x))\prod_{i:\alpha_i \in B}x^{h(B)}$$
\noindent and
$$C'x^{q(B',g')}=\prod_{i:\alpha'_i\not\in B'}\ini (\alpha'_j(x)-\alpha'_i(x))\prod_{i:\alpha'_i \in B'}x^{h(B')},$$
\noindent where $\alpha_j $ is a fixed element of  $B$.
Hence  $C=C'$, 
$q(B,g)=q(B',g')$ and $F_{B,g}(z)=F_{B',g'}(z+a_B)$.
By Lemma~\ref{weighted} the initial Newton polynomial of the discriminant depends only on the
critical values of $F_{B}(z)$ and on $q(B)$ for $B$ from the tree model. This proves Theorem~\ref{tc3} 
in transverse case. 

\medskip
\noindent To prove Theorem~\ref{tc3} in the case when $\ell_1=0$ and $\ell_2=0$ are tangent it is enough 
to take a smooth curve $\ell_3=0$ which is transverse to $\ell_1\ell_2f=0$ and apply the previously 
proved part to pairs of morphisms $(\ell_1,f)$, $(\ell_3,f)$ and $(\ell_3,f)$, $(\ell_2,f)$.
\end{proof}

\begin{Example}
\label{ejemJ}
Let $f=(y^2-x^2)^2+2x^4$. The discriminant
of $(x,f)$ is degenerate while the discriminant of $(x+y,f)$
is non-degenerate. The second discriminant can be easily
computed after the change of variables $x=x'-y'$, $y=y'$.
\end{Example}

\medskip
\noindent
{\small Evelia Rosa Garc\'{\i}a Barroso\\
Departamento de Matem\'atica Fundamental\\
Facultad de Matem\'aticas, Universidad de La Laguna\\
38271 La Laguna, Tenerife, Espa\~na\\
e-mail: ergarcia@ull.es}

\medskip

\noindent {\small Janusz Gwo\'zdziewicz,  Andrzej Lenarcik\\
Department of Mathematics\\
Technical University \\
Al. 1000 L PP7\\
25-314 Kielce, Poland\\
e-mails: matjg@tu.kielce.pl, ztpal@tu.kielce.pl}

\end{document}